\documentclass[12pt]{article}

\usepackage{amsfonts,amsmath, amssymb}

\usepackage{xcolor}

\usepackage{enumerate}
\usepackage{epsf,epsfig,amsfonts,a4wide}
\usepackage{amsfonts, mathdots}
\usepackage{amsmath,amssymb,amsthm}
\usepackage{url}
\usepackage{amsmath,amssymb,amsthm}

\usepackage{amssymb}
\usepackage{amsthm}
\usepackage{epsf,epsfig,amsfonts,graphicx}

\newtheorem{Theorem}{Theorem}[section]
\newtheorem{Lemma}{Lemma}[section]
\newtheorem{Proposition}{Proposition}[section]

\newtheorem{Ex}{Example}[section]
\newtheorem{Remark}{Remark}[section]

\theoremstyle{remark}

\newcommand{\be}{\begin{equation}}
\newcommand{\ee}{\end{equation}}

\newcommand{\R}{\mathbb{R}}\newcommand{\Id}{\textrm{\rm Id}}

\newcommand{\TT}{T}

\newcommand{\ddd}{\mathrm{d}}

\newcommand{\pd}[2]{\frac{\partial#1}{\partial#2}}

\newcommand{\dd}{{\mathrm d}\,}

\newcommand{\ii}[0]{\mathrm{i \,}}

\newcommand{\tr}{\operatorname{tr}}

\newcommand{\mathsfp}{{\mathsf p}}
\newcommand{\mathsfq}{{\mathsf q}}

\newcommand{\trace}{\operatorname{tr}}

\newcommand{\weg}[1]{}

\title{{ Applications of Nijenhuis geometry: Nondegenerate singular points of Poisson-Nijenhuis structures}}
\author{Alexey V. Bolsinov\footnote{ School of Mathematics,
 Loughborough University,
 LE11 3TU, UK   and   Faculty of Mechanics and Mathematics, Moscow State University, 119992, Moscow Russia\ \
 \quad {\tt A.Bolsinov@lboro.ac.uk} } \quad
\& \quad  Andrey Yu. Konyaev\footnote{Faculty of Mechanics and Mathematics, Moscow State University, 119992, Moscow Russia
 \ \ \quad {\tt  maodzund@yandex.ru}} \quad \& \quad Vladimir S. Matveev\footnote{
Institut f\"ur Mathematik, Friedrich Schiller Universit\"at Jena,
07737 Jena Germany  \ \ \quad {\tt  vladimir.matveev@uni-jena.de}} 
}  
\date{}

\begin{document}

\maketitle

\begin{abstract}
{We study and completely describe 
pairs of compatible Poisson structures near singular points of the recursion operator satisfying natural non-degeneracy condition.}
\end{abstract}



\section{Introduction}\label{sect:1}

This paper is a follow-up of our previous work
 \cite{Nijenhuis1} that can be considered as an introduction to Nijenhuis Geometry.  
Nijenhuis geometry studies 
 Nijenhuis operators, i.e., fields of endomorphisms $L=(L^i_j)$  of the tangent space 
with vanishing Nijenhuis torsion: that is, for any vector fields $\xi, \eta$
$$
N_L(\xi,\eta) = L^2[\xi,\eta] - L[L\xi,\eta] - L[\xi,L\eta] + [L\xi, L\eta]=0.   
$$
 The paper   \cite{Nijenhuis1} suggests  a programme  and  general strategy for studying such operators (and demonstrates that the programme  is realistic by  proving a series of non-trivial results). Important part of the strategy is the  study of  singular points of Nijenhuis operators. 

  The goal of this paper is to illustrate how methods, ideas  and results of  \cite{Nijenhuis1} work in  describing singular points  of more complicated geometric structures involving Nijenhuis operators as one of their ingredients   (cf. \cite[Sections 6.2 and 6.3]{Nijenhuis1})  and also  to give a proof for Theorem 6.3 announced in \cite{Nijenhuis1}.

The geometric structure we consider  in this paper came from 
  the theory of integrable systems, where Nijenhuis operators serve  as recursion operators for bihamiltonian structures \cite{gelzak3, kossman, magri1, magri2}, see also \cite{Fer1, Fer2}.

Recall that two  Poisson structures $P$ and $\widetilde P$   are called \emph{compatible}, if their sum $P + \widetilde P$ is also a Poisson structure. If $P$ and $\widetilde P$  are  
non-degenerate and therefore come from symplectic forms $\omega =P^{-1}$ and $\widetilde \omega=\widetilde P^{-1}$, then the compatibility condition  is equivalent to the property that the {\it recursion operator} $L$ given by the relation 
\begin{equation}
\label{omega12} 
\widetilde \omega(\,\cdot \,,  \cdot \,)=  \omega ( \,L \,\cdot \,,  \,\cdot  ) 
\end{equation}
is Nijenhuis.  Since $\widetilde \omega$ can be recovered from $\omega$ and $L$, we can reformulate the compatibility condition by saying that a symplectic structure $\omega$ and a Nijenhuis operator $L$ are {\it compatible} if

\begin{itemize}

\item[(I)]  the form $\widetilde \omega(\,\cdot \,,  \cdot \,)=  \omega ( \,L \,\cdot \,,  \,\cdot  )$ is skew-symmetric, i.e., is a differential 2-form,

\item[(II)] $\dd\widetilde\omega = 0$, i.e., this form is closed.

\end{itemize} 

In the case when $L$ is non-degenerate,  a compatible pair  $(\omega, L)$ defines a Poisson-Nijenhuis structure in the sense of  \cite{kossman, magri2}.  However, non-degeneracy of $L$ is not very essential as we can replace it with $L + c\cdot\Id$, $c\in\R$ and here we will think of compatible pairs  $(\omega, L)$ as a natural subclass of Poisson-Nijenhuis structures.

It is natural to ask about local simultaneous canonical form for compatible $\omega$ and $L$. In the case when $L$ is  algebraically generic at a point $\mathsfp\in M^{2n}$, i.e., its algebraic type (see the next section for details) remains the same in a certain neighbourhood $U(\mathsfp)$,  the answer was obtained in \cite{magri2}, \cite{turiel}, \cite{gelzak3}. We briefly recall some of these results.

Notice, first of all, that condition (I) imposes natural algebraic restrictions on the algebraic type of $L$: in the Jordan decomposition of $L$ all the blocks can be partitioned into pairs of equal blocks  (i.e. of the same size and with the same eigenvalue). In particular,  the characteristic polynomial of $L$ is a full square and each eigenvalue has even multiplicity.

If  $L$ has $n=\frac{1}{2}\dim M$ distinct real eigenvalues,  each of multiplicity 2,  and in addition their differentials are linearly independent at a point $p\in M$, then there exists a local symplectic coordinate system $x_1,\dots, x_n, \ p_1, \dots , p_n$ in which $L$ is diagonal with $x_i$ being its $i$-th eigenvalue \cite{magri2}:  
\begin{equation} 
\label{eq:GZform}
\omega = \sum_i dp_i\wedge dx_i, \ \ L= \operatorname{diag}(x_1,x_2, \dots ,x_n,\, x_1,x_2,\dots ,x_n).
\end{equation}
Equivalently, one can say that the pair  $(\omega, L)$ is the direct product of two-dimensional blocks $(\omega_i, L_i)$ of the form $\omega_i=dp_i\wedge dx_i$,  $L_i = x_i \cdot \Id$. For further use and without pretending that it is new, let us slightly generalise this result (cf. \cite{magri2}, \cite{turiel}, \cite{gelzak3}).  Consider 4 elementary examples of compatible pairs:
\begin{itemize}
\item[Type 1.]  $L$ has one real non-constant eigenvalue of multiplicity 2:
$$
\omega = \ddd p\wedge \ddd x, \quad L= \lambda(x,p) \cdot \Id.
$$

\item[Type 2.]  $L$ has one real constant eigenvalue of multiplicity $2k$;
$$
\omega = \sum_{j=1}^k \ddd p_j\wedge \ddd x_j, \quad L = \lambda \cdot \Id, \ \lambda\in\R. 
$$

\item[Type 3.] $L$ has one pair of non-constant complex conjugate eigenvalues  of multiplicity $2$:
$$
\omega = \mathrm{Re}  \, (\ddd z \wedge \ddd w), \quad L = \alpha(z,w) \cdot \Id + \beta(z,w) \cdot J,
$$
where $z= x+\ii y, w=u + \ii v$, $J$ denotes the corresponding complex structure and 
$\alpha(z,w) + \ii \beta(z,w)$ is a holomorphic function in $z$ and $w$.
\item[Type 4.]
$L$ has one pair of constant complex conjugate eigenvalues  of multiplicity $2k$:
$$
\omega = \mathrm{Re}  \, \left( \sum_{j=1}^{k}\ddd z_j \wedge \ddd w_j \right), \quad L = \alpha  \cdot \Id + \beta \cdot J,
$$
where $z_j= x_j+\ii y_j, w_j=u_j + \ii v_j$, $J$ denotes the corresponding complex structure and 
$\alpha+\ii\beta \in \mathbb C$, $\beta\ne 0$.
\end{itemize}

\begin{Theorem}[Folklore] \label{thm:gz0} 
Let $\omega$ and $L$ be compatible. Suppose that $L$ is semisimple and algebraically generic in a neighbourhood of $\mathsfp\in M$.  Then the pair $(\omega, L)$ locally splits into a direct product of `elementary blocks'  of 4 types described above.  If  $\ddd\lambda (\mathsfp) \ne 0$ or $\ddd(\alpha(\mathsfp) + \ii\beta (\mathsfp))\ne 0$ for some real or complex eigenvalue, then in the corresponding `elementary block'  we may set $\lambda = x$ (see Type 1) and $\alpha + \ii\beta = z$ (see Type 3).
\end{Theorem}

Notice that  to admit a compatible symplectic partner,  a semisimple algebraically generic Nijenhuis operator $L$ should satisfy    the following additional 
 condition:   non-constant eigenvalues of $L$ must be all of multiplicity two. 
 
If $L$ is algebraically generic but not necessarily semisimple, the description (rather non-trivial) of compatible pairs $(\omega, L)$  was obtained by Turiel \cite{turiel} under some additional assumptions on the differentials of $\tr L^k$, $k=1,\dots, n$.   These assumptions basically mean that each eigenvalue is either constant or its differential does not vanish.

The goal of the present paper is to give a simultaneous canonical form for $(\omega, L)$ at typical singular points, where the algebraic type of $L$ changes.

{\bf Acknowledgements.}  The work of A.\,Bolsinov and A.\,Konyaev was supported by Russian Science Foundation  (project 17-11-01303). 
Some results were obtained during the visits of A.\,Bolsinov and A.\ Konyaev to Jena supported by DAAD (via Ostpartnerschaft programm) and by the University of Jena, and during the research in pairs stay of A.\, Bolsinov, A.\, Konyaev and V. \, Matveev at  CIRM Trento. 
V.\, Matveev  thanks DFG for partial support via projects MA 2565/4
and MA 2565/6.

\section{Statement of the main result}\label{sect:2}

If we want to study singularities of Poisson-Nijenhuis structures (cf. Problem 5.17 in \cite{openprob}), then it is natural to ask which Nijenhuis operators admit compatible symplectic structures and what is a simultaneous canonical form for $\omega$ and $L$ near a (possibly singular) point $\mathsfp \in M$?

Let us explain in more detail what kind of singularities we have in mind.  For a symplectic structure $\omega$, all points $\mathsfp\in M$ are obviously equivalent to each other (in other words, there are no singular points). However it is not the case for an operator $L$ since
its algebraic properties may essentially vary from point to point.  More precisely,  at each point $\mathsfp\in M$,  we can define the {\it algebraic type}  (or {\it Segre characteristic}, see e.g. \cite{Fraser}) of  the operator $L(\mathsfp):  T_{\mathsfp} M \to T_{\mathsfp} M$  as the structure of its Jordan canonical form which is characterised by the sizes of Jordan blocks related to each eigenvalue $\lambda_i$ of $L(\mathsfp)$ (the specific values of $\lambda_i$'s are not important here). Following \cite[Definition 2.7]{Nijenhuis1}, we call a point $\mathsfp\in M$ {\it algebraically generic}, if the algebraic type of $L$ does not change in some neighbourhood $U(\mathsfp)\subset M$.  Otherwise, we say that $\mathsfp$ is {\it singular}.  

 Study of singularities of Nijenhuis and other related structures is an important part of the programme suggested in 
\cite{Nijenhuis1} (see also \cite{konyaev}). Actually, in most branches of geometry understanding  of singular points of relevant geometric structures plays a central role.

Among singular points,  we distinguish (see \cite[Definition 2.10]{Nijenhuis1})
an important subclass of the so-called {\it differentially non-degenerate} singular points $\mathsfp\in M$ which are characterised by the property that the differentials of the functions $\tr L^k$, $k=1,\dots, \dim M$, are linearly independent at $p$.   Instead of  the traces of powers  of $L$ one can equivalently consider the coefficients $\sigma_1, \dots, \sigma_m$, $m=\dim M$, of the characteristic polynomial of $L$:
$$
\chi_L(t)= \det (t\cdot \Id - L) = t^m + \sigma_1 t^{m-1} + \dots + \sigma_{m-1} t + \sigma_m. 
$$
Relation (I), however,  implies that the eigenvalues of  $L$ have always  even multiplicity, which in turn implies that $L$ cannot
 be differentially non-degenerate (since the characteristic polynomial is a full square and has at most $\frac{1}{2} \dim M$ functionally independent  coefficients). 

As a natural analog of differential non-degeneracy in this setup  we suggest the condition that the differentials of the functions $\trace {L}$,  $\trace {L^2}$,\dots, 
$\trace {L^n}$, $n=\frac{1}{2}\dim M$,  are linearly independent at a given (singular) point $\mathsfp\in M$.  In terms of the characteristic polynomial of $L$ which is a full square { (of a polynomial with smooth coefficients)} in our case:
$$
\chi_L(t)=(t^n + h_1 t^{n-1} + \dots + h_{n-1} t + h_n)^2,   \qquad m=2n=\dim M,
$$
this condition amounts to linear independence of the differentials $\ddd h_1(\mathsfp), \dots, \ddd h_n(\mathsfp)$.
Our goal is to study local canonical forms for $\omega$ and $L$ at such points. This problem  is clearly important for finite dimensional integrable systems and we do hope that it might also be important for understanding of  bihamiltonian structures in infinite dimension appearing in the theory of integrable ODEs and PDEs.  The main result of this paper is the following theorem whose weaker version is 
announced in \cite{Nijenhuis1}.

\begin{Theorem} \label{thm:gz} Let $\omega$ and $L$ be compatible  (i.e., define a Poisson-Nijenhuis structure on $M^{2n}$) and real analytic. Suppose that  at a point $\mathsfp\in M^{2n}$,  the differentials  $\ddd\tr L$, \dots, $\ddd\tr L^n$ are linearly independent.
Then  there exists a local coordinate system $x_1,\dots,x_n, p_1,\dots,p_n$ such that $\omega=\sum_i \ddd x_i\wedge \ddd p_i$ and $L$ is given by the matrix 
\begin{equation} \label{AS1} \begin{pmatrix} A & 0_n\\ S &  A^\top\end{pmatrix},\end{equation}
where 
	\begin{equation} \label{AS2}
	A =\begin{pmatrix}  -x_1&    1   &   0   &   \cdots      &  0       \\
                 -x_2&  0     &     1 &             \ddots   &   \vdots     \\  
						  \vdots & \vdots &\ddots & \ddots                &  0  \\
                -x_{n-1}& \vdots &   & \ddots &1 \\
								-x_n& 0      &\cdots &   \cdots        &            0       
\end{pmatrix}, \ \ \   
S=\begin{pmatrix}  
								  0    & \!\! -p_2    &\!\! -p_3 &   \cdots &   \!\!    -p_n       \\
								  p_2  & 0       & 0 &     \cdots &      0        \\  
                  p_3  &  0 & \vdots &             &   \vdots  \\ 
								  \vdots&  \vdots &  \vdots       &      &  \vdots     \\  
							    p_n  &   0       &0 &	 \cdots 	  &    0       
\end{pmatrix},      \end{equation} 
$0_n$ is the zero $n\times n$-matrix, and $A^\top$ denotes the transposed of $A$.
\end{Theorem}

\begin{Remark}{\rm
We would like to point out that this theorem mainly concerns singular points of Poisson-Nijenhuis structures, i.e., those at which the algebraic type of $L$ changes  (recall that the case of  algebraically generic points has basically been covered by the results of 
 \cite{magri2}, \cite{turiel}, \cite{gelzak3}).  The typical and most important example of such a  situation is a ``perturbation of a Jordan block''  when all the eigenvalues of $L$  simultaneously vanish at the point $\mathsfp\in M^{2n}$.  In this situation, the $x$-coordinates of $\mathsfp$  are $(0,0,\dots, 0)$. However, the formulas \eqref{AS1} and \eqref{AS2} are universal in the sense that they are applicable to both algebraically generic and singular points of  any algebraic kind under the only condition  that the differentials  $\ddd\tr L$, \dots, $\ddd\tr L^n$ are linearly independent.

 More precisely, 
 if a Nijenhuis-Poisson structure $(L,\omega)$ satisfies the conditions of Theorem \ref{thm:gz} and in addition, at the point $\mathsfp\in M^{2n}$, the operator $L$ has $s$ distinct eigenvalues $\lambda_1,\dots,\lambda_s$  (some of which can be complex conjugate) of multiplicities $2k_1,\dots,2k_s$,  then the local canonical form for $L$, $\omega$ is given by \eqref{AS1} and \eqref{AS2}  where the $x$-coordinates $x_1^0,\dots, x^0_n$ of the point $\mathsfp$  are the coefficients of the polynomial
 $$
 \prod_{i=1}^s  (t-\lambda_i)^{k_i} = t^n + x_1^0 t^{n-1} +x_2^0 t^{n-2} + \dots + x_{n-1}^0 t + x_n^0,
 $$
that is, the square root of the characteristic polynomials $\chi_{L(\mathsfp)}$.
 }\end{Remark}

The rest of the paper is devoted to the proof of Theorem \ref{thm:gz}.  First we discuss one cohomological property of $\mathrm{gl}$-regular operators serving as a key step in our proof and also interesting on its own.  This property uses some existence and uniqueness results from general theory of PDEs that we clarify in Appendix A.  The proof of Theorem \ref{thm:gz} is given 
in Section \ref{sect:4}. In Appendix B we formulate an alternative version of Theorem \ref{thm:gz}.

All the objects we are dealing with (functions, manifolds, maps) are assumed to be real analytic.

\section{One cohomology-like property of $\mathrm{gl}$-regular Nijenhuis operators.}\label{sect:3}

Let $V$ be a finite-dimensional vector space. We will say that an operator $A: V \to V$ is $\mathrm{gl}$-regular, if the dimension of its conjugacy class $\mathcal O(A)=\{ PAP^{-1}~:~ P\in {\mathrm{GL}}(V)\}$ is maximal (i.e.,  equals  $n^2-n$). Equivalently, this regularity condition means that in an appropriate basis the matrix of $A$ takes the so-called companion form as in \eqref{eq:companion}.  Typical examples of $\mathrm{gl}$-regular operators include a semisimple operator with distinct eigenvalues and a single Jordan block. 

 To clarify the relation between the $\mathrm{gl}$-regularity condition and our main Theorem  \ref{thm:gz}, it is worth mentioning that the linear independence of the differentials $\ddd\tr A(\mathsfp), \dots, \ddd\tr A^n(\mathsfp)$  automatically implies that $A(\mathsfp)$ is 
$\mathrm{gl}$-regular\footnote{This is a purely algebraic fact, a particular case of the famous Kostant theorem \cite[Theorem 0.1]{Kostant}, which can be easily verified in our case related to the Lie algebra $\mathrm{gl}(n)$.}. Notice that $A$ and $L$ from Theorem \ref{thm:gz} are related by  \eqref{AS1},  hence $\tr L^k = 2\tr A^k$ so that the assumptions of this theorem imply $\mathrm{gl}$-regularity  of the block $A$. These two operators $A$ and $L$ are related via the Turiel extension construction discussed in Appendix B.  However in this section, the reader should think of $A$ just as an operator defined on an $n$-dimensional manifold $Q$ and not related to $L$ in any sense.

\begin{Proposition}
\label{prop:cohomol}
Let $A$ be a (real-analytic) Nijenhuis operator which is $\mathrm{gl}$-regular at a point $\mathsfp \in Q$.
Assume that $\Omega$ is a closed (real-analytic) $2$-form such that the form 
\begin{equation}
\label{eq:num1}
\Omega_A(\xi, \eta) = \Omega(A\xi, \eta) + \Omega(\xi, A\eta)
\end{equation}
is also closed.
Then locally in some neighbourhood of $\mathsfp \in Q$ there exists a (real-analytic) function $U$ such that 
\begin{equation}
\label{eq:num2}
\Omega = \dd (A^* \dd U).
\end{equation}
\end{Proposition}

\begin{Remark}{\rm
Note that the $\mathrm{gl}$-regularity condition is essential. It can be easily checked that the statement of Proposition 
\ref{prop:cohomol} fails for the operator $A =\begin{pmatrix}  x_1 & 0 \\ 0 & x_2 \end{pmatrix}$ at those points $\mathsfp \in Q$ where the eigenvalues collide, i.e.,  $x_1 = x_2$. }\end{Remark}

\begin{Remark}{\rm 
Proposition \ref{prop:cohomol} can be interpreted in terms of the Nijenhuis differential \cite{kossman}, \cite{magri}. Denote by $\Lambda^k (Q)$ the module of differential $k-$forms on a manifold $Q$, $\Lambda^0$ are functions on $Q$ (definition works for both smooth and analytic cases) and $\ddd: \Lambda^k \to \Lambda^{k + 1}$ is standard differential of $k-$forms.

Given { an operator $A$}, one can define a new differential $\ddd_A: \Lambda^k \to \Lambda^{k + 1}$ which is completely determined by the following properties (see  \cite[Appendix]{magri})
\begin{equation}\label{coh}
    \begin{aligned}
    &1. \ \ddd_Af = A^* \ddd f, \quad f \in \Lambda^0, \\
    &2. \ \ddd_A (\alpha + \beta) = \ddd_A \alpha + \ddd_A \beta, \quad \alpha, \beta \in \Lambda^k, \\
    &3. \ \ddd_A (\alpha \wedge \beta) = \ddd_A \alpha \wedge \beta + (-1)^k \alpha \wedge \ddd_A \beta, \quad \alpha \in \Lambda^k, \beta \in \Lambda^s , \\
    &4. \ \ddd_A \ddd + \ddd \ddd_A = 0.
    \end{aligned}
\end{equation}

If  $\mathcal N_A = 0$, we have that $\ddd_A^2 = 0$ (see  \cite[Corollary 6.1]{kossman} for $E = TQ$, $N = A$ and $\mu$ being standard commutator of vector fields). 

 The explicit formula for the operator  $\dd_A:\Lambda^k \to \Lambda^{k + 1}$ is as follows \cite{kossman}: 
$$
\dd_A  = [i_A, \dd],
$$
where $i_A: \Lambda^m \to \Lambda^{m}$  (for any $m\in \mathbb N$) is defined by 
$$
(i_A  \alpha) ( \xi_1, \xi_2, \dots, \xi_m) =\alpha (A\xi_1, \xi_2, \dots, \xi_m) +
\alpha (\xi_1, A \xi_2, \dots, \xi_m) + \dots + \alpha (\xi_1, \xi_2, \dots, A\xi_m),
$$
so that for 2-forms $i_A$ is essentially given by \eqref{eq:num1}, i.e.  $i_A \Omega = \Omega_A$.    We assume that $i_A(f)=0$ for any function $f$. 
}\end{Remark}

Thus,  Proposition \ref{prop:cohomol} can be reformulated as follows: given $2-$form $\Omega$ such that $\ddd \Omega = 0$ and $\ddd_A \Omega = 0$ (the latter formula follows from  (\ref{eq:num1}) combined with \eqref{coh}), one can always find a function $U$ such that $\Omega = \ddd \ddd_A U$.

{ The goal of the  rest of  this section is to prove Proposition \ref{prop:cohomol}.}  
Recall that we are mainly interested in those points $\mathsfp\in Q$ where $A$ changes its algebraic type.  The key example of such a situation is  a perturbation of a Jordan block, i.e.,  $A(\mathsfp)$ is similar to a single Jordan block  with zero eigenvalue, but generically the operator $A(\mathsfq)$,  $\mathsfq\in V(\mathsfp)$ is semisimple with distinct eigenvalues (as in \eqref{AS2}). 

We however first consider the case when $A(\mathsfp)$ is regular and semisimple. 

\begin{Lemma}\label{lem:num1}
In the assumptions of Proposition \ref{prop:cohomol}, let $A(\mathsfp)$ be semisimple with distinct eigenvalues. Then the statement of Proposition \ref{prop:cohomol} holds. Moreover, the general solution $U$ of \eqref{eq:num2} depends on (is parametrised by) $n$ arbitrary functions of one variable.
\end{Lemma}

\begin{proof} Consider a local coordinate system $y_1,\dots, y_n$ in which $A$ takes diagonal form
$$
A = \operatorname{diag}\bigl(\lambda_1(y_1), \dots, \lambda_n(y_n)\bigr).
$$

Since $A$ is Nijenhuis,  such a coordinate system exists  \cite{haant, nij}.

Consider a closed 2-form $\Omega = \sum_{i<j} \omega_{ij} (y) \ddd y_i \wedge \ddd y_j$  satisfying \eqref{eq:num1},  i.e.  
$\Omega_A =\sum_{i<j} (\lambda_i + \lambda_j)\omega_{ij} (y) \ddd y_i \wedge \ddd y_j$ is also closed.

We need to solve \eqref{eq:num2}  which takes, in coordinates $y_1,\dots, y_n$, the following form:
\begin{equation}
\label{eq:num4}
(\lambda_j - \lambda_i)  \frac{\partial ^2 U}{\partial y_i y_j}   =   \omega_{ij},
\end{equation}
or, equivalently,  
$$
 \frac{\partial ^2 U}{\partial y_i  \partial y_j}   =   \frac{ \omega_{ij}}{\lambda_j -\lambda_i}, \qquad  i\ne j.
$$

This equation is consistent if and only if
$$
\frac{\partial}{\partial y_k} \left(\frac{ \omega_{ij}}{\lambda_i -\lambda_j}\right)  = \frac{\partial}{\partial y_i} \left(\frac{ \omega_{jk}}{\lambda_j -\lambda_k}\right) =
\frac{\partial}{\partial y_j} \left(\frac{ \omega_{ki}}{\lambda_k -\lambda_i}\right) 
$$
for $i\ne j \ne k \ne i$,  which can also be rewritten as
$$
\frac{\partial \omega_{ij}}{\partial y_k} \left(\lambda_i -\lambda_j\right)^{-1}  = 
\frac{\partial \omega_{jk}}{\partial y_i} \left(\lambda_j -\lambda_k\right)^{-1} =
\frac{\partial \omega_{ki}}{\partial y_j} \left(\lambda_k -\lambda_i\right)^{-1} 
$$

These two linear relations follow from (and in fact, are equivalent to)
$$
\frac{\partial \omega_{ij}}{\partial y_k} +
\frac{\partial \omega_{jk}}{\partial y_i} +
\frac{\partial \omega_{ki}}{\partial y_j} =0  \quad \mbox{(closedness of $\Omega$)}
$$
and
$$
\frac{\partial \omega_{ij}}{\partial y_k} (\lambda_i+\lambda_j) +
\frac{\partial \omega_{jk}}{\partial y_i} (\lambda_j +\lambda_k) +
\frac{\partial \omega_{ki}}{\partial y_j} (\lambda_k+ \lambda_i)=0 \quad \mbox{(closedness of $\Omega_A$)}.
$$

Hence, \eqref{eq:num4}  (or equivalently \eqref{eq:num2}),  admits a solution $U_0$.  Since \eqref{eq:num4} is linear in $U$,  its general solution is defined modulo solutions of the corresponding homogeneous system    
$$ 
\frac{\partial ^2 U}{\partial y_i \partial y_j}=0,
$$ 
that is,
 $$
 U = U_0(y_1,\dots,y_n) + u_1(y_1) + \dots + u_n(y_n).
 $$ 
In other words, $U$ depends on $n$ functions in one variables, as stated.
\end{proof}

\begin{Remark}\label{rem:11}{\rm
On can similarly check that  Proposition \ref{prop:cohomol} holds under the additional condition that $\mathsfp\in Q$ is algebraically generic, i.e. the algebraic type of $A$ remains locally constant in a neighbourhood of $\mathsfp\in Q$.  
The splitting theorem \cite[Section 3.2]{Nijenhuis1} allows one to assume without loss of generality that $A$ has either a single real eigenvalue or single pair of complex conjugate eigenvalues.  In this case, due to $\mathrm{gl}$-regularity,  $A$ will be similar to a single Jordan block (real or complex)  and the general solution of \eqref{eq:num2} can be found explicitly by using the canonical forms for such operators from \cite{Nijenhuis1}.  We omit details of this computation here, although the explicit description of functions $U$ satisfying the condition $\dd (A^* \dd U)=0$ for $A$ being a Nijenhuis Jordan block, could be interesting for several reasons. We will provide such a description in another paper.  
}\end{Remark}

The following statement is a simple fact from Linear Algebra.  It is well known that every $\mathrm{gl}$-regular  operator $A$ can be reduced  (at one single point $\mathsfp\in Q$)
to companion form\footnote{In fact, for $\mathrm{gl}$-regular Nijenhuis operators $A$, this can be done simultaneously  for all points from some  neighbourhood of $p$ by choosing an appropriate local coordinate system \cite{Nijenhuis3}.}, i.e., 
\begin{equation}
\label{eq:companion}
A(\mathsfp) = A_{\mathrm{comp}} =\begin{pmatrix} 
-\sigma_1 & 1 & & \\
\vdots &0 & \ddots & \\
-\sigma_{n-1} & & \ddots & 1 \\
-\sigma_n & 0 & \dots & 0
\end{pmatrix}.
\end{equation}
where $\sigma_i$ will be automatically the coefficients of the characteristic polynomial $\chi_A(t) = \det (t\cdot \Id - A)$.

\begin{Lemma} 
\label{lem:num2}
Let $x_1,\dots, x_n$ be a local coordinate system in a neighbourhood of $p\in Q$ such that in these coordinates $A(p) =  A_{\mathrm{comp}}$.  Then 
PDE system \eqref{eq:num2} can be resolved w.r.t. all second order derivatives of $U$ except those of the form $U_{x_i x_n}$ ($i=1,\dots, n$), i.e., can be rewritten in the form
\begin{equation}\label{eq:num3}
U_{x_i x_j} := \frac{\partial ^2 U}{\partial x_i \partial x_j} =  h_{ij} (x_1,\dots,x_n,  U_{x_1}, \dots, U_{x_n},  U_{x_1x_n}, \dots, U_{x_n x_n}), \quad i,j=1,\dots, n-1.
\end{equation}
\end{Lemma}

\begin{proof}   At the point $p\in Q$ itself,  the proof is a simple exercise in Linear Algebra. In all neighbouring points, representation \eqref{eq:num3} for system  \eqref{eq:num2} holds  by continuity. \end{proof}

In order to   prove Proposition \ref{prop:cohomol}  we will show that PDE system \eqref{eq:num3} possesses sufficiently many solutions  (namely, the space of solutions is parametrised by $n$ functions in one variable).  To that end, we use the following key property of PDE systems.

\begin{Lemma}\label{lem:num3}
The following properties of the PDE system  \eqref{eq:num3} are equivalent:
\begin{itemize}

\item[{\sf (i)}]  for any initial straight line $x_1=a_1, \dots, x_{n-1}=a_{n-1}$, $a_i\in\R$ and any initial conditions given on it
$$
\begin{aligned}
U(a_1,\dots, a_{n-1}, x_n) &= v(x_n),\\
U_{x_1}(a_1,\dots, a_{n-1}, x_n) &= v_1(x_n),\\
\qquad \vdots \\
U_{x_{n-1}}(a_1,\dots, a_{n-1}, x_n) &= v_{n-1}(x_n),\\
\end{aligned}
$$
there locally exists a unique real analytic solution $U(x_1,\dots, x_n)$ of  \eqref{eq:num3}.

\item[{\sf (ii)}] The following compatibility conditions identically hold in virtue of \eqref{eq:num3}:
\begin{equation}
\label{eq:num6}
D_{x_k} h_{ij} = D_{x_i} h_{jk} = D_{x_j}  h_{ki}.
\end{equation}
\end{itemize}
\end{Lemma}

\begin{proof}
The proof will be given in Appendix A. 
\end{proof}

In our case, the functions $h_{ij}$  are linear (but not homogeneous) in the derivatives $U_{x_i}$ and $U_{x_ix_n}$, i.e., 
$$
h_{ij} = \sum_{s=1}^n A_{ij}^s(x)  U_{x_s x_n} +  \sum_{s=1}^n B_{ij}^s(x) U_{x_s}   + C_{ij}(x)
$$

For the given operator $A$ and $2$-form $\Omega$,  the compatibility conditions \eqref{eq:num6} take the form
\begin{equation}
\label{eq:num10}
c_\alpha (x) \equiv 0, \qquad \alpha\in I,
\end{equation}
where $c_\alpha$ are some functions in coordinates $(x_1,\dots, x_n)$ that can be explicitly written in terms of the components of the operator $A$ and  2-form $\Omega$ and their partial derivatives, i.e.,
$$
c_\alpha(x) = \tilde c_\alpha \left(  A_i^j, \, \frac{\partial A_i^j}{\partial x_k},  \, \omega_{ij}, \, \frac{\partial \omega_{ij}}{\partial x_k}     \right)
$$

 Recall that both $A$ and $\Omega$ are rather special and satisfy very strong conditions. However,  straightforward verification of compatibility conditions \eqref{eq:num10}  (and even explicit description of them in terms of $A$, $\Omega$ and their derivatives)
 is not an easy task.  Instead, we will use the fact that $c_\alpha (x)$ is  real-analytic in $x$ and for this reason it is sufficient to verify these conditions only for some open subset  $V_0\subset V(\mathsfp)$.  As such a subset we consider a small neighbourhood $U(\mathsfq)\subset V(\mathsfp)$ of a point $q\in V(p)$ at which $A(\mathsfq)$ is semisimple (similar verification can be done for $q$ being algebraically generic\footnote{Here we use the fact that algebraically generic points form an open and everywhere dense subset.}, cf. Remark \ref{rem:11}).  According to Lemma \ref{lem:num1}, the general solution $U$ of PDE system \eqref{eq:num2}  (which is equivalent to  \eqref{eq:num3})  depends on $n$ arbitrary real analytic functions of one variable.  Basically, this property is equivalent to condition {\sf (i)}  from Lemma \ref{lem:num3} and therefore the compatibility conditions \eqref{eq:num6} from {\sf (ii)}  hold.

However, the collections of $n$ functions parametrising the space of solutions of \eqref{eq:num2}  in Lemma \ref{lem:num1} and in Lemma \ref{lem:num3} are formally different.  So we need to justify the fact that by choosing appropriate functions $u_1(y_1), \dots, u_n(y_n)$  (see Lemma   \ref{lem:num1})  we can fulfil arbitrary initial conditions 
$v(x_n), v_1(x_n), \dots, v_{n-1}(x_n)$  from Lemma \ref{lem:num3}.

Consider the parametrised straight line $\gamma$ given in coordinates $x_1, \dots, x_n$  (from Lemma \ref{lem:num3}) as 
$$
\gamma(t):  \quad  x_1 =a_1, \ x_2 = a_2, \ \dots \ , \ x_{n-1}=a_{n-1}, \ x_n = t.
$$
The same line in coordinates $y_1, \dots, y_n$ (from Lemma \ref{lem:num1}) will be given as
$$
\gamma(t):  \quad  y_1 =\phi_1(t), \ y_2 = \phi_2(t), \ \dots \ , \ y_{n-1}=\phi_{n-1}(t), \ y_n = \phi_n(t).
$$
It can be easily seen that each of these functions $\phi_i$ is locally invertible, i.e., $\dfrac{d\phi_i}{dt}\ne 0$ at $t=0$. Indeed, this property means that the entries of the last column $\xi$ of the Jacobi matrix $J=\Bigl(  \frac{\partial y}{\partial x} \Bigr)$ (at the origin) are all different from zero.  This matrix satisfies the relation $A_{\mathrm{comp}}(x) = J^{-1} A_{\mathrm{diag}} (y)J$ or, equivalently, $J A_{\mathrm{comp}} = A_{\mathrm{diag}} J$.  The latter can be interpreted as follows:  the columns of $J$ are (in the reverse order) of the form 
$$
\xi, \ A_{\mathrm{diag}} \xi, \  A_{\mathrm{diag}}^2 \xi, \ \dots \  A_{\mathrm{diag}}^{n-1} \xi.
$$    
In particular, these vectors are linearly independent.  For a diagonal matrix $A_{\mathrm{diag}}$ this condition holds if and only if all the components of $\xi\in\R^n$ are different from zero, as stated. 
 
The initial conditions {\sf (i)}  from Lemma \ref{lem:num3}  can be rewritten in the form
$$
U(\gamma(t)) = v(t), \ U_{x_1}(\gamma(t)) = v_1(t), \ \dots \ , \ U_{x_{n-1}}(\gamma(t)) = v_{n-1}(t). 
$$

It is more convenient to rewrite the first condition in the form $U_{x_n}(\gamma(t)) = v_n(t)= v'(t)$, $U(\gamma(0)) = c$ (where $c\in \R$ is an arbitrary constant). In other words, we replace $v(t)$ with its derivative $v'(t)$ and the value $c\in\R$ at $t=0$.  

Using the standard relations $U_{x_j}=\sum_\alpha \frac{\partial y_\alpha}{\partial x_j} U_{y_\alpha}$ we can rewrite these conditions in coordinates $y_1,\dots, y_n$ as follows:
$$
\left.\sum_\alpha \frac{\partial y_\alpha}{\partial x_j}  U_{y_\alpha}\right|_{\gamma(t)} = v_j(t),
$$
and using  $U(y_1,\dots, y_n)  = U_0(y_1,\dots, y_n) + u_1(y_1) + \dots + u_n(y_n)$: 
$$
\left.\sum_\alpha \frac{\partial y_\alpha}{\partial x_j}  \bigl((U_0)_{y_\alpha}  + u_\alpha' \bigr)\right|_{\gamma(t)}= v_j(t)
$$
or
$$
u_\alpha'  (\phi_\alpha(t)) =  b_\alpha(t) + \sum_\beta   a^j_\alpha (t) v_j(t), \quad \mbox{where } \
b_\alpha(t) =  - (U_0)_{y_\alpha} (\gamma(t)), \  a^j_\alpha (t)=\frac{\partial x_j}{\partial y_\alpha}(\gamma(t)). 
$$

Taking into account the fact  that $\phi_\alpha(t)$ is locally invertible, we see that the derivatives of the functions $u_1, \dots, u_n$ can be uniquely reconstructed from $v_1, \dots, v_n$ and vice versa as required. 

We have shown that in a neighbourhood of any point $\mathsfq$ at which $A(\mathsfq)$ is diagonalisable over $\R$,   we can find a solution of \eqref{eq:num2} with arbitrary initial conditions {\sf (i)}  from Lemma \ref{lem:num3}. This guarantees the fulfilment of the stong compatibility conditions \eqref{eq:num6} from {\sf (ii)} in this neighbourhood and therefore (due to real analyticity) in the whole neighbourhood $V(\mathsfp)$.   The latter, in turn,  guarantees the existence of local solutions of \eqref{eq:num2} in a neighbourhood of $\mathsfp\in Q$, completing the proof of Proposition \ref{prop:cohomol}.

\section{Proof of the main theorem}\label{sect:4}

Recall that we consider $\omega$ and $L$ to be compatible and real analytic. We also assume that at a point $\mathsfp\in M^{2n}$  all the eigenvalues of $L$ vanish,  but the differentials  $\ddd\tr L$, \dots, $\ddd\tr L^n$ are linearly independent.
    
Since  the bilinear form 
$\widetilde\omega(\cdot,\cdot) = \omega( L\cdot,\cdot)$ is skew-symmetric, the characteristic polynomial of $L$ is a full square:  
\begin{equation}
\label{eq:fullsquare}
\chi_L(t)= (t^n + h_1 t^{n-1} + \ldots  + h_n)^2. 
\end{equation}
i.e., the square of a certain polynomial of degree $n$. The condition that  
   the differentials of the functions $\trace {L}$, \dots ,  $\trace {L^n}$ are linearly independent implies that the differentials of the functions   $h_1, \dots,h_n$ are linearly independent too and we can take them as the first $n$ coordinates $x_1, \dots,x_n$ 
	of a local coordinate system. It is well known that these functions commute  with respect to the Poisson bracket related to $\omega$ (see \cite{magri2, gelzak3}) which, by  Darboux theorem, implies local  existence of  functions $p_1,\dots,p_n$ such that $(x_1, \dots , x_n, p_1, \dots , p_n)$ is a canonical coordinate system for $\omega$.

In these coordinates, the operator $L$ takes an ``almost required'' form.  To see this we use the following formula proved in \cite{Nijenhuis1}.   Let $\sigma_1, \dots, \sigma_{m}$, $m=\dim M$, be the coefficients of the characteristic polynomial of $L$:
$$
\chi_L(t) = t^{m} + \sigma_{1} t^{m-1} + \sigma_{2} t^{m-2}+ \dots  + \sigma_{m-1} t + \sigma_m   
$$	
considered as smooth functions on $M$. Then in any coordinate system $y_1, \dots, y_{m}$, the matrix $L(y)$ satisfies the following relation 
$$
J L(y) = S J,  \quad\mbox{where } S=\begin{pmatrix}   - \sigma_1 & 1 & & \\ \vdots & 0 &\ddots & \\ -\sigma_{m-1} &\vdots &\ddots & 1 \\ -\sigma_m & 0 &\dots & 0 \end{pmatrix} 
$$	
and $J = \left(  \dfrac{\partial \sigma_i}{\partial y_j} \right)$ is the Jacobi matrix of the functions $\sigma_1,\dots, \sigma_m$ w.r.t.  $y_1,\dots, y_m$.	In our case,  the characteristic polynomial $\chi_L(t)$  is a full square as in \eqref{eq:fullsquare} so that its coefficients $\sigma_i$ ($i=1,\dots, 2n$) are some explicit polynomials in $h_\alpha=x_\alpha$ ($\alpha=1,\dots, n$) and do not depend on the other half of coordinates $p_\alpha$ ($\alpha=1,\dots, n$).  It is now straightforward to check that in coordinates $(x,p)$  
the first $n$ rows of the operator $L$ are as in  (\ref{AS1},\ref{AS2}). Furthermore, the algebraic compatibility condition  (I) implies that the lower right $n\times n$-block of $L$ is transposed to the upper left $n\times n$-block, so that $L$ takes the form
\begin{equation} 
\label{eq:L} 
L= \begin{pmatrix} 
A & 0 \\ 
\widehat S & A^\top
\end{pmatrix}
\end{equation}  
with $A$ as in (\ref{AS2}) and $\widehat S$ being just a  skew-symmetric matrix (again due to condition (I)) whose components  may {\it a priori}   depend on all variables. 
 
 Next we use the differential compatibility condition (II) saying that the 2-form $\widetilde\omega(\cdot,\cdot) = \omega(L\cdot,\cdot)$ is closed.  This form  is given by     
 \begin{equation}
\label{calc}
\widetilde\omega= \sum_{i=1}^n x_i \, \ddd x_1\wedge \ddd p_i + \sum_{i,j=1}^n \widehat S_{ij} \, \ddd x_i \wedge \ddd x_j
\end{equation} 
and its differential 
is 
\begin{equation}
\label{calc1}
\dd\widetilde\omega =  - \sum_{i=2}^n \ddd x_1\wedge \ddd x_i \wedge \ddd p_i + \sum_{i,j,k=1}^n \frac{\partial \widehat S_{ij}}{\partial x_k} \, \ddd x_k \wedge \ddd x_i \wedge \ddd x_j +
\sum_{i,j,k=1}^n \frac{\partial \widehat S_{ij}}{\partial p_k}\, \ddd p_k \wedge \ddd x_i \wedge \ddd x_j . 
\end{equation}
Substituting  $\widehat  S =   S  +  \TT$ in the formula above, where $S$ is  as in (\ref{AS2}), we observe that the first sum  in (\ref{calc1})
cancels with the terms coming from $S$  and we obtain     
\begin{equation}\label{calc2}
\dd \widetilde\omega=\sum_{i,j,k=1}^n\frac{\partial \TT_{ij}}{\partial p_k} \, \ddd p_k \wedge \ddd x_i \wedge \ddd x_j \
 + \sum_{i,j,k=1}^n \frac{\partial \TT_{ij}}{\partial x_k} \, \ddd x_k \wedge \ddd x_i \wedge \ddd x_j  .
\end{equation}
Since $\widetilde\omega$ is closed, we see that $\TT_{ij}$ does  not depend  on the $p$-variables, so that the form $\sum_{i,j=1}^n\TT_{ij}(x) dx_i\wedge  dx_j$  can be viewed as a closed 2-form $T$ on the local $n$-dimensional coordinate chart $(x_1,\dots,x_n)$. 

To complete the proof,  we will show that the form $\TT$ can be ``killed'' by a suitable canonical transformation of the form 
\begin{equation}
\label{eq:coord}
(X_1,\dots,X_n, P_1,\dots,P_n)=  \left(x_1,\dots, x_n, \ p_1+ \tfrac{\partial U}{\partial x_1},\dots, p_n+ \tfrac{\partial U}{\partial x_n}\right),\end{equation} 
where $U$ is a   function of $x_1,\dots,x_n$ (called sometimes {\it generating function}).  This transformation preserves the  (canonical) symplectic structure $\omega$ so that we only need to look after the change of $L$.

The Jacobi matrix of coordinate transformation \eqref{eq:coord} is
$
J=\begin{pmatrix} \mathrm{id}& 0\\ \ddd^2 U & \mathrm{id}\end{pmatrix}
$,
 where $\mathrm{id}$ denotes the identity $n\times n$ matrix and $\ddd^2U= \Bigl(\tfrac{\partial^2 U}{\partial x_i\partial x_j}\Bigr)$ is the Hessian matrix of $U$. Hence, after  transformation \eqref{eq:coord} the matrix of $L$ given by \eqref{eq:L} takes the form  
 $$ 
 L_{\textrm{new}}= JLJ^{-1}=  \begin{pmatrix} \mathrm{id}& 0\\ \ddd^2U & \mathrm{id}\end{pmatrix} \begin{pmatrix} A & 0 \\ \widehat S & A^\top\end{pmatrix}\begin{pmatrix} {\mathrm{id}}& 0\\ -\ddd^2U & \mathrm{id}\end{pmatrix}= \begin{pmatrix} A & 0 \\ \widetilde  S  & A^\top\end{pmatrix}, 
 $$
where $\widetilde S = \widehat S +  \ddd^2U \cdot A - A^\top \!\! \cdot \ddd^2U$. Taking into account the transformation of the $p$-coordinates in $S$, we see that $L_{\textrm{new}}$ takes the required form (\ref{AS1},\ref{AS2}) if and only if the generating function $U$ satisfies the following equation:
\begin{equation}\label{eq:U}
0      =   \TT + \ddd^2U \cdot A - A^\top \!\! \cdot \ddd^2U + \dd U \wedge \ddd x_1,
\end{equation}  
 which can be  equivalently rewritten in a more invariant form: 
 \begin{equation}\label{eq:U1}
   \dd(A^*\dd U)=\TT 
\end{equation}
that coincides with the PDE system \eqref{eq:num2} treated in the previous section.  In order to check the existence of a generating function $U$ solving \eqref{eq:U1}, it remains to verify two conditions from Proposition \ref{prop:cohomol} imposed on the 2-form $\TT$.

They both follow from the fact that  
$L= \begin{pmatrix} A & 0 \\ \widehat S & A^\top\end{pmatrix}$ with $\widehat S = S + T$ and $A$ and $S$ as in \eqref{AS2}, is a Nijenhuis operator. It is easily checked by a straightforward computation that vanishing Nijenhuis torsion of $L$  is exactly equivalent to the two geometric conditions we need, namely  that the  2-forms  
 \begin{equation} \label{eq:forms}
T =  \sum_{i,j=1}^n\TT_{ij} \ddd x_i \wedge \ddd x_j \ \  \mbox{and}  \ \ T_A=\sum_{i,j,k=1}^n \left({A^k}_i \TT_{kj}+ 
   {A^k}_j \TT_{ik}\right)\ddd x_i \wedge \ddd x_j\end{equation}
     are closed.  This allows us to apply Proposition \ref{prop:cohomol} that now guarantees the existence of  $U$ solving \eqref{eq:U1}  and hence completes the proof of Theorem \ref{thm:gz}.

\section{Appendix A: Proof of Lemma \ref{lem:num3}.}

The goal in this Appendix is to prove Lemma \ref{lem:num3} concerning the solvability of  the PDE system
$$
U_{x_i x_j}  =  h_{ij} (x_1, \dots, x_n,  U_{x_1}, \dots, U_{x_n},  U_{x_1x_n}, \dots, U_{x_n x_n}), \quad i,j=1,\dots, n-1.
$$

Lemma \ref{lem:num3} is not a new fact and can be understood as a special case of the Cartan-K\"ahler theorem
(see e.g. \cite{Gold}). We suppose that  people working  in the Cartan-K\"ahler theory would find it trivial. However,  deriving Lemma \ref{lem:num3} from a very general and  non-trivial Cartan-K\"ahler theorem seems to require more work than proving it  directly.  Unfortunately, we have not found this lemma in the required form in the literature, so we will prove it here in order to make the paper self-contained.

As one usually does in the theory of differential equations, we set $f_1 = U_{x_1}, f_2 = U_{x_2},  \dots , f_n = U_{x_n}$ and rewrite the above system in the form 
\begin{equation}
\label{eq:num11}
F_{x_i} = H_{i}(x, F, F_{x_n}), \quad \mbox{with } F = \begin{pmatrix}  f_1 \\ \vdots \\ f_k \end{pmatrix},  \quad i=1,\dots, n-1.
\end{equation}

In our case the number of unknown functions  $f_j$ equals the number of variables $x_i$,  i.e., $k=n$ but in general this assumption is irrelevant. 

In this view,  we will prove an analogue of Lemma \ref{lem:num3}  for a more general PDE system \eqref{eq:num11} (from which Lemma \ref{lem:num3} can be easily derived by using the above substitution $f_i = U_{x_i}$). Recall that all the functions involved are assumed
to be real analytic. 

\begin{Lemma}\label{lem:5.1} 
The following properties of \eqref{eq:num11} are equivalent:

\begin{itemize}

\item  for any initial straight line $x_1=a_1, \dots, x_n = a_n$ and any real analytic initial conditions on it 
$$
\begin{aligned}
f_1(a_1,\dots, a_{n-1} , x_n) &= v_1(x_n)\\ 
f_2(a_1,\dots, a_{n-1} , x_n) &= v_2(x_n)\\
\dots &\\
f_k(a_1,\dots, a_{n-1} , x_n) &= v_k(x_n)\\
\end{aligned}
$$
there locally exists a unique real analytic solution $F = \begin{pmatrix}  f_1 \\ \vdots \\ f_n \end{pmatrix}$ of \eqref{eq:num11}.
\item The following compatibility conditions are identically fulfilled in virtue of  \eqref{eq:num11}:
$$
D_{x_i} H_j = D_{x_j} H_i, \quad i,j=1,\dots, n-1.
$$ 

\end{itemize}
\end{Lemma}

\begin{proof}

We give a proof in the simplest non-trivial case\footnote{The general case is not essentially different.}, namely,   for three variables $x, y, z$ and one function $f$ so that  \eqref{eq:num11} becomes
\begin{equation}
\label{eq:num12}
f_{x} = h_1(x,y,z, f, f_z) \quad \mbox{and} \quad f_{y} = h_2(x,y,z, f, f_z).
\end{equation}
 
 Assume that  the condition $D_x h_2 = D_y h_1$ holds identically in virtue of the system.  In more details:
$$
D_x h_2  = \frac{\partial h_2}{\partial x} + \frac{\partial h_2}{\partial f} f_x +  \frac{\partial h_2}{\partial f_z} f_{zx} =
\frac{\partial h_2}{\partial x} + \frac{\partial h_2}{\partial f} h_1 +  \frac{\partial h_2}{\partial f_z} D_zh_1(x,y,z, f, f_z) =
$$

$$
\frac{\partial h_2}{\partial x} + \frac{\partial h_2}{\partial f} h_1 + \frac{\partial h_2}{\partial f_z} \left(
\frac{\partial h_1}{\partial z} + \frac{\partial h_1}{\partial f} f_z +  \frac{\partial h_1}{\partial f_{z}} f_{zz}\right) =
$$

And similarly 
$$
D_y h_1  = \frac{\partial h_1}{\partial y} + \frac{\partial h_1}{\partial f} h_2 +  \frac{\partial h_1}{\partial f_z} \left(
 \frac{\partial h_2}{\partial z} + \frac{\partial h_2}{\partial f} f_z +  \frac{\partial h_2}{\partial f_z} f_{zz}\right)
$$

In these two expressions we see an additional formal variable $f_{zz}$.  In this simples case the corresponding coefficients in the both take the form    $\frac{\partial h_1}{\partial f_z}  \frac{\partial h_2}{\partial f_z}$ and therefore coincide automatically. However in general it is not necessarily the case so that the equality these coefficients should be considered as an extra condition additional to 

$$
\frac{\partial h_2}{\partial x} + \frac{\partial h_2}{\partial f} h_1 +  \frac{\partial h_2}{\partial f_z} \left(
 \frac{\partial h_1}{\partial z} + \frac{\partial h_1}{\partial f} f_z  \right) = \frac{\partial h_1}{\partial y} + \frac{\partial h_1}{\partial f} h_2 +  \frac{\partial h_1}{\partial f_z} \left(
 \frac{\partial h_2}{\partial z} + \frac{\partial h_2}{\partial f} f_z   \right)
 $$

The identical fulfilment of this condition means that the functions in the right hand side and the left hand side coincide as functions of 5 independent variables $x,y,z, f$ and $f_z$. 

Without loss of generality, we will assume that the parameters $a_1$ and $a_2$ in the initial condition vanish. For given initial condition $f(0,0,z) = v(z)$ we construct the desired solution in two steps.  We first construct a unique function of two variables $u(y,z)$  that is a local real-analytic solution of the Cauchy problem
$$
u_y = h_2(0,y,z, u, u_z), \quad u(0,z) = v(z).
$$
The existence of $u(y,z)$ follows from the classical Cauchy-Kovalevskaya theorem. If the desired solution $f$ of \eqref{eq:num12} exists then necessarily $f(x,0,z)=u(x,z)$. 

The next step is solving the Cauchy problem
$$
f_x = h_1(x,y,z,f,f_z), \quad f(0,y,z)=u(y,z).
$$
Again, a unique local real-analytic solution $f(x,y,z)$ exists by Cauchy-Kovalevskaya theorem.

Thus we have constructed a function that satisfies the initial condition and solves the first equation of \eqref{eq:num12}, as for the the second equation,  it is fulfilled only on the plane $x=0$. Our goal is to show that $D_x h_2 = D_y h_2$ guarantees that the second equation holds at all the points. 

\begin{Lemma}
Let $f$ satisfy the conditions
$$
f_x (x,y,z) = h_1(x,y,z,f,f_z)   \quad \mbox{and} \quad  f_y (0,y,z) = h_2(0,y,z,f,f_z).
$$
Then if the compatibility condition $D_x h_2 = D_y h_1$ is  identically  fulfilled, then  
$$
f_y (x,y,z) = h_2(x,y,z,f,f_z) \quad\mbox{for all } x,y,z.
$$
\end{Lemma}

\begin{proof}
Consider the functions  
$$
\widehat h_2(x,y,z) = h_2(x,y,z,f(x,y,z),f_z(x,y,z)) \quad\mbox{and}\quad  \widehat h_1(x,y,z) = h_1(x,y,z,f(x,y,z),f_z(x,y,z))
$$ 
and compute the partial derivative of $f_y - \widehat h_2$ w.r.t. $x$:
$$
\frac{\partial}{\partial x} \left( f_y - \widehat h_2\right) 
= \frac{\partial }{\partial y} ( f_x ) -    \frac{\partial h_2}{\partial x}  -  \frac{\partial h_2}{\partial f}  f_x  -  \frac{\partial h_2}{\partial f_z}  f_{zx} =
$$
$$
\frac{\partial \widehat h_1}{\partial y} -   \frac{\partial h_2}{\partial x}  -  \frac{\partial h_2}{\partial f} \, \widehat h_1  -  \frac{\partial h_2}{\partial f_z}  
\frac{\partial \widehat h_1}{\partial z} =
$$
$$
\frac{\partial h_1}{\partial y}  + \frac{\partial h_1}{\partial f}  f_y + \frac{\partial h_1}{\partial f_z}  (f_y)_z  -   \frac{\partial h_2}{\partial x}  -  \frac{\partial h_2}{\partial f} \, \widehat h_1  -  \frac{\partial h_2}{\partial f_z}  
\frac{\partial \widehat h_1}{\partial z} =
$$
and taking into account the compatibility condition:
$$
 \frac{\partial h_1}{\partial f}  f_y + \frac{\partial h_1}{\partial f_z}  (f_y)_z    -  \frac{\partial h_1}{\partial f} \, \widehat h_2  -  \frac{\partial h_1}{\partial f_z}  
\frac{\partial \widehat h_2}{\partial z} =
 \frac{\partial h_1}{\partial f}  \left(f_y-\widehat h_2\right) + \frac{\partial h_1}{\partial f_z} \left(f_y-\widehat h_2\right)_z.
$$

We see that the function $f_y - \widehat h_2$ satisfies the linear PDE:
$$
\frac{\partial}{\partial x} \left( f_y - \widehat h_2\right) = \frac{\partial h_1}{\partial f}  \left(f_y-\widehat h_2\right) + \frac{\partial h_1}{\partial f_z} \left(f_y-\widehat h_2\right)_z.
$$
Moreover for $x=0$  we have $f_y(0,y,z) - \widehat h_2(0,y,z) =0$.  This implies that $f_y -\widehat h_2$ vanishes identically, that is,
$$
f_y(x,y,z) = h_2(x,y,z, f(x,y,z), f_z(x,y,z)),
$$
as stated.
\end{proof}

The proof of the converse statement is obvious.  Indeed, if $f$ is the solution of \eqref{eq:num12} with  initial condition $f(a_1, a_2, z) = u(z)$, then  $f_{xy} = f_{yx}$ implies that
$D_y H_1 (x,y,z, f, f_z) = D_x H_2 (x,y,z, f, f_z)$.  Moreover, since  $u(z)$ and $u_z(z)$  can be arbitrarily chosen for fixed $x=a_1, y=a_2$ and $z$,  we see that  $D_y H_1 = D_x H_2$ holds identically  (with $x$, $y$, $z$, $f$ and $f_z$ being formal independent variables), as required.  \end{proof}

\section{Appendix B: Turiel extension construction and an alternative form
of the main theorem.}

Formulas (\ref{AS1},\ref{AS2})  from Theorem  \ref{thm:gz} can be naturally interpreted in terms of the Turiel extension construction 
 \cite{turiel1}. This gives an additional proof that they indeed define a Nijenhuis operator that is compatible with the canonical symplectic  structure $\sum_i \ddd x_i\wedge \ddd p_i$.

Recall that for an arbitrary Nijenhuis operator $A$ on a manifold $Q$,  the Turiel construction provides a natural extension of $A$ to the cotangent bundle  $M^{2n}=T^*Q^n$ as follows. Consider  local coordinates $x_1, \dots, x_n$ on $Q$  and denote by $x_1, \dots, x_n, p_1, \dots, p_n$ the corresponding canonical coordinates on $T^*Q$. Next, consider the operator on $T^*Q$  given in these coordinates by 
\begin{equation}\label{ext}
    L = \left(\begin{array}{cc}
         A &  0 \\
         S & A^\top
    \end{array}\right),
\end{equation}
where 
\begin{equation} \label{tur:ext} 
S_{ij} := \sum \limits_{\alpha = 1}^n p_{\alpha} \left(\pd{A^{\alpha}_{i}}{x_j} - \pd{A^{\alpha}_{j}}{x_i} \right).
\end{equation}

By \cite{turiel1}, this operator is  Nijenhuis and compatible with the canonical symplectic structure $\sum_i \ddd p_i\wedge \ddd x_i$ on $T^*Q$. 

\begin{Ex} {\rm
By direct calculations we see that  the ``normal form'' from  Theorem \ref{thm:gz} is in fact  
 the Turiel extension of the Nijenhuis operator in the ({\it first}) companion form (\ref{eq:companion}) with $\sigma_i= x_i$: indeed, $S$  and $A$ from \eqref{AS2} are related by \eqref{tur:ext}. 
}\end{Ex} 

Observe  that the Turiel construction is evidently geometric and does not depend on the choice of a coordinate system $x_1, \dots, x_n$ in which we apply it.  Indeed, $L$ is the recursion operators for two 2-forms $\omega$ and $\widetilde\omega$ each of which is geometric (see \cite{IbortMagri}):  the first of them is canonical $\omega = - \ddd \theta$, $\theta = \sum p_\alpha \ddd x_\alpha$, the second is defined by the following invariant formula $\widetilde\omega = - \ddd\widetilde\theta$ with $\widetilde\theta=
A^* \theta = \sum A^i_j (x) p_i\ddd x_j$.
One therefore may look for coordinate systems on $Q$ such that in the corresponding coordinate system on $T^*Q$ 
 the extended operator $L$ takes a more convenient form. The interpretation  of  ``convenient'' may depend on the problem we are solving and 
on the situation in which the recursion  Nijenhuis operator appeared.  As an illustration, let us show that one can choose a coordinate system on $Q$ 
in such a way that $S$ in \eqref{ext} vanishes.

 We say that a Nijenhuis operator $A$ is in \emph{the second companion form} if its matrix, in local coordinates $y_1, \dots, y_n$,  has the form   (cf. \eqref{eq:companion})
\begin{equation}\label{second}
    A = \left(
\begin{array}{ccccc}
     0 & 1 & 0 & \dots & 0  \\
     0 & 0 & 1 & \dots & 0 \\
     \vdots & \vdots & \ddots & \ddots & \vdots \\
     0 & 0 &\dots & 0 & 1 \\
     - \sigma_n & - \sigma_{n-1} &  \dots & \dots & - \sigma_1 \\
\end{array}
\right).
\end{equation}

 Not every operator of the form (\ref{second}) is Nijenhuis.  The next lemma gives necessary and sufficient condition for this property to hold.

\begin{Lemma}\label{lem1} An operator $A$  given by (\ref{second}) is  Nijenhuis, if and only if 
\begin{eqnarray} 
    \ddd(A^*\ddd y_n) &= & 0 \label{sig1}\\
	\ddd ({A^*}^2\ddd y_n) &=&0. \label{sig2}
\end{eqnarray}
\end{Lemma}

Condition (\ref{sig1}) simply means that the differential 1--form given by the last row of $A$ is closed, i.e., $\frac{\partial  \sigma_{n-i+1}}{\partial  y_j} =  
\frac{\partial  \sigma_{n-j+1}}{\partial  y_i}$
 for any $i,j$;  it is a linear condition on $\sigma_i$.   The other condition (\ref{sig2})  is nonlinear.

\begin{proof}   We first recall (see e.g. \cite[Def. 2.5]{Nijenhuis1})
that for an operator $A$, its  Nijenhuis tensor $\mathcal N_A$ viewed as a mapping from 1--forms to  2--forms  is given by  \begin{equation}
\label{eq:NijForForms}
\mathcal N_A: \ \alpha \mapsto  \beta( \cdot\, , \cdot ):= \ddd({A^*}^2\alpha) (\cdot\, ,\cdot) + \ddd\alpha (A\cdot\, ,A\cdot) -
\ddd(A^*\alpha) (A\cdot\, , \cdot) - \ddd(A^*\alpha) (\cdot\, , A\cdot). 
\end{equation}
If $A$ is Nijenhuis and  both $\alpha$ and $A^*\alpha$  are closed, this  formula implies that also $ {A^*}^2\alpha$ is closed. 
For our $A$ given by (\ref{second}),   the 
 forms $\ddd y_1$ and $\ddd y_2=A^*\ddd y_1$ are closed. Thus, if $A$ is Nijenhuis, also the form $({A^*})^{n}\ddd y_1=A^*\ddd
y_n$ is closed, which gives us (\ref{sig1}). Similarly, the form   $({A^*})^{n+1}\ddd y_1={A^*}^2\ddd y_n$ is closed, which gives us (\ref{sig2}). 

In order to prove Lemma \ref{lem1} in the other direction, let us observe that since $\mathcal N_A$ is a tensor,  it is sufficient to check that (\ref{eq:NijForForms}) vanishes for the basis forms $\alpha = \ddd y_i$, $i=1,\dots, n$.  For them  each term in (\ref{eq:NijForForms}) is zero.  \end{proof}

\begin{Ex} \label{ex:noS}
{\rm The Turiel extension of the Nijenhuis operator $A$ in the second companion form is 
\begin{equation}\label{noS}
    L = \left(\begin{array}{cc}
         A &  0_n \\
         0_n & A^\top
    \end{array}\right),
		\end{equation}
where
$0_n$ is the zero $n \times n$ matrix. Indeed, for fixed $i$ the 1-form  ${A^{\alpha}_{i}}\ddd x_\alpha$ is just ${A^*}^{(i-1)}\ddd x_1$; so the corresponding sum in 
(\ref{tur:ext}) vanishes for $i<n$ trivially and for $i=n$
by (\ref{sig1}).  
}\end{Ex}

We use this observation to give another canonical form for Poisson-Nijenhuis structures from Theorem \ref{thm:gz}. To that end we use a coordinate transformation  $x \mapsto y$ that reduces the Nijenhuis operator $A$ from Theorem \ref{thm:gz}   (being in the {\it first} companion form):
\begin{equation} \label{AS}
	A =\begin{pmatrix}  -x_1&    1   &   0   &   \cdots      &  0       \\
                 -x_2&  0     &     1 &             \ddots   &   \vdots     \\  
						  \vdots & \vdots &\ddots & \ddots                &  0  \\
                -x_{n-1}& \vdots &   & \ddots &1 \\
								-x_n& 0      &\cdots &   \cdots        &            0       
\end{pmatrix}. \end{equation} 
to the {\it second} companion form \eqref{second}. 
Suitable new coordinates  $(y_1, \dots ,y_n)$ are easy to find by using \cite[Proposition 2.1]{Nijenhuis1} which states 
that for any Nijenhuis operator $A$ and any $k\in \mathbb{N}$ we have 
\begin{equation}\label{eq:fromN1}
\frac{1}{k}\,\ddd \tr (A^k)= \frac{1}{k-1}\,A^*\ddd \tr (A^{k-1}).
\end{equation}
Let us set:
$$ 
y_1= -\tr  A, \ \ y_2=-\frac{1}{2}\tr  
 A^2, \  \dots , \  y_n=-\frac{1}{n}\tr  A^n.
$$ 
Notice that the transition formulas $x \mapsto y$ express the traces $y_1, \dots, y_n$ of powers of $A$ (with some coefficients)  in terms of the coefficients of its characteristic polynomial $x_1, \dots, x_n$. Such a transformation is known to be invertible   and polynomial in both directions (see below the formulas for the inverse transformation $y \mapsto x$).

From \eqref{eq:fromN1} we see that $A^*\ddd y_i= \ddd y_{i+1}$ 
(for $i\in 1,\dots,n-1$) meaning that the first $n-1$ rows of $A$ in the $y$-coordinate system  are as we want, i.e., 
 in these coordinates $A$ takes the second companion form \eqref{second}. The corresponding functions $\sigma_k(y)$  in the last row of \eqref{second}   can be  explicitly constructed. Indeed, $\sigma_k$ are the coefficients of the characteristic polynomial   $\chi_A(t)=\det (t\cdot \operatorname{Id} - A ) = t^n+ \sigma_1 t^{n-1}+\dots+\sigma_n$  so that the  functions $\sigma_k(y_1,\dots, y_n)$ coincides with $x_k$ and are  given by the condition 
\begin{equation} 
\label{G} 
\sigma_k\left(-\tr  A,-\frac{1}{2}\tr  A^2,\dots,-\frac{1}{n}\tr  A^n\right)=
\textrm{coefficient of $t^{n-k}$ in   $\chi_A(t)$.}
\end{equation}      
This condition is well known and 
 the functions  $\sigma_i$  satisfying \eqref{G} can be obtained from the famous  Newton-Girard polynomials (by appropriate rescaling and changing sings).  They depend neither on the matrix $A$  nor on the dimension  and are given by the following   universal formula, see e.g. \cite{wiki}:  
   \begin{equation} \label{eq:wiki}\sigma_{k}=\sum_{\begin{array}{c}{m_{1}+2m_{2}+\cdots +km_{k}=k}\\ {m_{1}\geq 0,\ldots ,m_{k}\geq 0} \end{array}}\prod _{i=1}^{k}{\frac {y_{i}^{m_{i}}}{m_{i}!{}}}\end{equation} 

For example, the first five functions $\sigma_1,\dots,\sigma_5$   are given by
\begin{eqnarray*}\sigma_1&= & y_1\\
{ \sigma_2} & = &y_{{2}}+
\frac{1}{2}\,{y_{{1}}}^{2} \\ 
{\sigma_3}& =&y_{{3}}+y_{{1}}y_{{2}}+\frac{1}{6}\,{y_{{1}}}^{3}\\ 
{ \sigma_4}  &= & y_{{4}}+\frac{1}{24}\,{y_{{1}}}^{4}+\frac{1}{2}\,y_{{2}}{y_{{1}}}^{2}+y_{{1}}y_{{3}}
+\frac{1}{2}\,{y_{{2}}}^{2}\\
 \sigma_5 & =&y_{{5}}+{\frac{1}{120}{{y_{{1}}}^{5}}}+\frac{1}{6}\,y_{{2}}{y_{{1}}}^{3}+\frac{1}{2}\,
y_{{3}}{y_{{1}}}^{2}+\frac{1}{2}
y_{{1}}{y_{{2}}}^{2}+y_{{1}}y_{{4}}+y_{{3}}y
_{{2}}
\end{eqnarray*}

Now, if we naturally extend the coordinate transformation $x\mapsto y$  from $Q$ to the cotangent bundle $M=T^*Q$,   then in new coordinates $(y, p_y)$ the operator $L$ will be the Turiel extension of $A$ in the second companion form  \eqref{second}, whereas $\omega$  will remain canonical.  Taking into account Example \ref{ex:noS}, we get the following alternative version of  Theorem \ref{thm:gz}.

 \begin{Theorem}  Let $L$ and $\omega$ be compatible (i.e., define a Poisson-Nijenhuis structure on $M^{2n}$) and real analytic. Suppose that at a point $p \in M^{2n}$  the differentials $\ddd \operatorname{tr} L, \dots, \ddd \operatorname{tr} L^n$ are linearly independent. Then there exists a local coordinate system $y_1, \dots, y_n, p_1, \dots, p_n$ such that $\omega = \sum_i \ddd y_i \wedge \ddd p_i$ and $L$ is given by the matrix 
\begin{equation} \label{LL}
     L = \left(\begin{array}{cc}
         A &  0_n \\
         0_n & A^\top
    \end{array}\right),
\end{equation}
where $A$ is given by \eqref{second} with $\sigma_k$ given  by \eqref{eq:wiki}.
\end{Theorem}


\end{document}